\documentclass{amsart}
\usepackage{amsfonts,latexsym,amsmath,amscd,amssymb}
\usepackage{bbm}

\theoremstyle{plain}
\newtheorem{theorem}{Theorem}
\numberwithin{theorem}{section}

\numberwithin{corollary}{section}

\newtheorem{definition}{Definition}
\numberwithin{definition}{section}

\newtheorem{lemma}{Lemma}
\numberwithin{lemma}{section}

\newtheorem{proposition}{Proposition}
\numberwithin{proposition}{section}

\newtheorem{remark}{Remark}
\numberwithin{remark}{section}

\numberwithin{equation}{section}

\newcommand {\be}{\begin{equation}}
\newcommand {\ee}{\end{equation}}

\newcommand{\h}{\begin{eqnarray*}}
 \newcommand{\e}{\end{eqnarray*}}
\newcommand{\CC}{\mathbb{C}}

\begin{document}
\title[Elliptic Genera and Transgression]{Elliptic Genera, Transgression and Loop Space Chern-Simons Forms}

\author{Qingtao Chen}
\address{Q. Chen, Department of Mathematics, University of California,
Berkeley, CA, 94720-3840} \email{chenqtao@math.berkeley.edu}
\date{12 Jan, 2007}

\author{Fei Han}
\address{F. Han, \ Department of Mathematics, University of California,
Berkeley, CA, 94720-3840} \email{feihan@math.berkeley.edu}


\maketitle

\begin{abstract} We compute the Chern-Simons transgressed forms of
some modularly invariant characteristic forms, which are related to
the elliptic genera. We study the modularity properties of these
secondary characteristic forms and the relations among them. We also
compute the Chern-Simons forms of some vector bundles over free loop
space.

\end{abstract}

\section {Introduction}

Connections in vector bundles play a very important role in
differential geometry. The famous Chern-Weil theory relates
connections to the theory of characteristic classes, which provides
a geometric way to understand characteristic classes and therefore
play a fundamental role in global differential geometry. Let
$\nabla$ be a connection on an $n$-dimensional real or complex
vector bundle $E\rightarrow M$, where $M$ is a compact smooth
manifold. The basic local invariant of $\nabla$ is its curvature,
which is a closed 2-form with values in the $so(n)$-bundle or
$su(n)$-bundle associated to $E$. Invariant polynomials of the
curvature give us Chern-Weil characteristic forms, which are closed
and thus represent cohomology classes in the de Rham cohomology of
$M$. These classes don't depend on the choice of the connections and
are actually topological invariants of $E$. Particular invariant
polynomials give us the famous Chern classes and Pontrajin classes.

The Chern-Smions theory studies the dependence of characteristic
forms on $\nabla$, which, by transgression method, leads to
secondary geometric invariants, called the Chern-Simons forms. Chern
and Simons [10] are led to this theory by concrete geometric
questions in combinatorial and conformal geometry. Cheeger and
Simons go further along this clue to define some refined secondary
invariants called the Cheeger-Simons differential characters [9]. It
turns out that these secondary invariants are very useful in many
areas of mathematics and physics. For example, Witten [27] uses the
secondary invariant associated to a particular characteristic form
to construct a topological quantum field theory in 3 dimension and
obtains some quantum invariants including the Jones polynomial of
knots as well as new invariants of 3-manifolds. The Chern-Simons
forms or more generally the geometric invariants of connections have
special advantages to study flat vector bundles. A flat bundle is
the vector bundle equipped with a connection whose curvature
vanishes identically. Hence, by the Chern-Weil theory, all
characteristic forms vanish. Thus we cannot read off any information
from this theory. However based on the idea of Chern-Simons and the
method of transgression, one is able to construct certain cohomology
classes by using flat connections, which turn out to be useful tools
to study flat vector bundles. In Section 2 of this paper, we briefly
review the method of the Chern-Simons transgression and give Theorem
2.2, which serves as a convenient tool to do transgressions on the
characteristic forms that we are going to use.

A lesson we have learned from the above stories is that one may
apply the Chern-Simons transgression to particularly chosen
characteristic forms and obtain interesting secondary characteristic
forms, which might have potential applications. In this paper, we
apply the Chern-Simons transgression to some modularly invariant
characteristic forms and obtain modularly invariant secondary
characteristic characteristic forms. The characteristic forms that
we use here are related to the theory of elliptic genera. Ochanine
[24] first introduced the notion of elliptic genera from the
topological point of view. He was motivated by a study of
$S^1$-action on Spin manifolds by Landweber and Stong [19], who
attempted to answer a question raised by Witten [26] on rigidity of
certain twisted Dirac operators on Spin manifold. In [28] Witten
reinterpreted the Landweber-Stong elliptic genus as the index of the
formal signature operator on loop space as well as introduced the
formal index of the Dirac operator on loop space (called the Witten
operator), which is the Witten genus. In Section 3 of this paper, we
apply the Chern-Weil theory to express the local indices of some
formal elliptic operators on loop spaces by connections (or
curvatures), which are characteristic forms with modular invariance
properties. Then in Section 4, we apply the Chern-Simons
transgression to these characteristic forms to obtain some
interesting secondary characteristic forms. Most of these secondary
characteristic forms have modularly invariant properties except for
the one associated to the Witten operator. However it turns out that
this specific secondary characteristic form is a modular form over
$SL(2, \mathbb{Z})$ when it's associated to two flat connections on
a flat manifold. Moreover this secondary characteristic form cannot
be of weight 2 over $SL(2, \mathbb{Z})$, which agrees with the
standard fact in the theory of modular forms. As the original
characteristic forms, the secondary characteristic forms we obtain
are also modularly related. We would like to point out that the
modularities of these secondary characteristic forms are not a
direct consequence of the modularities of the original
characteristic forms. See details in Section 4. We hope that these
new geometric invariants of connections with modularity properties
obtained here could be applied somewhere.

Motivated by string theory, people have been attempting to
generalize many things like vector bundles, Dirac operators, the
Atiyah-Singer index theory and so on to loop spaces. The relevant
vector bundles over loop space are constructed from ordinary
finite-rank vector bundles $V\to M$. For example, in Section 5, we
will consider vector bundles $\mathcal{V}$ and $\mathcal{V}'$ over
loop space:
$$\mathcal{V}:=\bigotimes_{j=1}^{\infty}\Lambda_{-q^{j-{1/2}}}(V_\CC),\
\ \
\mathcal{V}':=\bigotimes_{j=1}^{\infty}\Lambda_{q^{j-{1/2}}}(V_\CC).$$
As Witten remarked in his lecture notes [29] that physically, this
analogue is very important because it arises in heterotic string
theory. Following [21], in our discussion section (Section 5), we
describe the Atiyah-Singer index theory on loop spaces. Then for
flat bundles on loop space, we compute their Chern-Simons forms,
which turn out to have modular properties. This could be understood
as certain generalization of the Chern-Simons forms for flat bundles
on finite dimensional smooth manifolds to free loop spaces.

\section{Construction of the Chern-Simons transgressed forms}
Let's briefly review the construction of the Chern-Simons
transgressed forms in this section.

Let $M$ be a compact smooth manifold and $T^*M$ denote the
cotagent bundle of $M$. We denote by $\Lambda^*(T^*M)$ the
(complex) exterior algebra bundle of $T^*M$, and $$\Omega^*(M,
\mathbb{C})\triangleq\Gamma(\Lambda^*(T^*M))$$ the space of smooth
sections of $\Lambda^*(T^*M)$. In particular, for any integer $p$
such that $0\leq p \leq \mathrm{dim}M$, we denote by
$$\Omega^p(M, \mathbb{C})\triangleq\Gamma(\Lambda^p(T^*M))$$
the space of smooth $p$-forms over $M$. Let $E$ be a smooth
complex vector bundle over $M$. We denote by $\Omega^*(M; E)$ the
space of smooth sections of the tensor product vector bundle
$\Lambda^*(T^*M)\otimes E$ obtained from $\Lambda^*(T^*M)$ and
$E$,
$$ \Omega^*(M; E)\triangleq\Gamma(\Lambda^*(T^*M)\otimes E).$$
Let $\mathrm{End}(E)$ be the bundle of endomorphisms of $E$. On $
\Omega^*(M; \mathrm{End}(E))$, one can define a super Lie bracket
(cf. [30]) by extending the Lie bracket operation on
$\mathrm{End}(E)$ as follows: if $\omega, \eta \in \Omega^*(M)$ and
$A, B\in \Gamma(\mathrm{End}(E))$, then we use the convention that
$$[\omega A, \eta B]=(\omega A)(\eta
B)-(-1)^{(\mathrm{deg}\omega)(\mathrm{deg}\eta)}(\eta B)(\omega
A).$$ It's not hard to see that: for any $A, B \in \Omega^*(M;
\mathrm{End}(E))$, the trace of $[A, B]$ vanishes (cf. [31, Lemma
1.7]).

Let $\nabla^{E}$ be a connection on $E$ and $A\in \Omega^*(M;
\mathrm{End}(E))$. One has the following result,

\begin{lemma}[cf. {[31, Lemma
1.8] }]The following identity holds, \be d\,
\mathrm{tr}[A]=\mathrm{tr}\left[[\nabla^E, A]\right].\ee
\end{lemma}

Let $$ f(x)=a_0+a_1x+\cdots+a_nx^n+\cdots$$ be a power series in one
variable. Let $R^E=(\nabla^E)^2$ be the curvature of $\nabla^E$ on
$E$. The trace of $f(R^E)$ is an element in $\Omega^*(M,
\mathbb{C})$. A form of the Chern-Weil theorem (cf. [31, Theorem
1.9]) can be stated as follows.

\begin{theorem} 1) The form $\mathrm{tr}[f(R^E)]$ is closed. That is,
$$d\mathrm{tr}[f(R^E)]=0;$$

2) If $\widetilde{\nabla}^E$ is another connection on $E$ and
$\widetilde{R^E}$ is its curvature. For any $t\in [0,1]$, let
$\nabla_t^E$ be the deformed connection on $E$ given by
$\nabla_t^E=(1-t)\nabla^E+t\widetilde{\nabla}^E.$ Let $R^E_t, t \in
[0,1]$, denote the curvature of $\nabla_t^E$. Let $f'(t)$ be the
power series obtained from the derivative of $f(x)$ with respect to
$x$. Then the following identity holds, \be
\mathrm{tr}[f(\widetilde{R}^E)]-\mathrm{tr}[f(R^E)]=d\
\int_0^1\mathrm{tr}\left[\frac{d\nabla^E_t}{dt}f'(R^E_t)\right]dt.\ee
\end{theorem}
\begin{proof} 1) From Lemma 2.1, one verifies directly that
$$ d
\mathrm{tr}[f(R^E)]=\mathrm{tr}\left[[\nabla^E, f(R^E)]\right]$$
$$=\mathrm{tr}\left[a_1[\nabla^E, R^E]+\cdots+a_n[\nabla^E,
(R^E)^n]+\cdots+\right]=0,$$ as for any integer $k\geq 0$ one has
the obvious Bianchi identity \be [\nabla^E, (R^E)^k]=[\nabla^E,
(\nabla^E)^{2k}]=0.\ee

2) Note that $\nabla_t^E$ is a connection on $E$ such that
$\nabla_0^E=\nabla^E$ and $\nabla_1^E=\widetilde{\nabla}^E.$
Moreover, $$ \frac{d\nabla_t^E}{dt}=\widetilde{\nabla}^E-\nabla^E
\in \Omega^1(M, \mathrm{End}(E)).$$

We deduce that
$$\frac{d}{dt}\mathrm{tr}\left[f(R^E_t)\right]=\mathrm{tr}\left[\frac{dR^E_t}{dt}f'(R^E_t)\right]=
\mathrm{tr}\left[\frac{d(\nabla^E_t)^2}{dt}f'(R^E_t)\right]$$
$$=\mathrm{tr}\left[\left[\nabla^E_t, \frac{d \nabla^E_t}{dt}\right]f'(R^E_t)\right]=
\mathrm{tr}\left[\left[\nabla^E_t, \frac{d
\nabla^E_t}{dt}f'(R^E_t)\right]\right],$$ where the last equality
follows from the Bianchi identity (2.3).

Combining with Lemma 2.1, we have \be
\frac{d}{dt}\mathrm{tr}\left[f(R^E_t)\right]=d\
\mathrm{tr}\left[\frac{d \nabla^E_t}{dt}f'(R^E_t) \right], \ee
from which one obtains
$$\mathrm{tr}[f(\widetilde{R}^E)]-\mathrm{tr}[f(R^E)]=d\
\int_0^1\mathrm{tr}\left[\frac{d\nabla^E_t}{dt}f'(R^E_t)\right]dt.$$

\end{proof}

The transgressed term \be \int_0^1\mathrm{tr}\left[\frac{d
\nabla^E_t}{dt}f'(R^E_t) \right]dt \ee is usually called a {\bf
Chern-Simons term}.

In particular, let $M$ be a 3-dimentional oriented compact smooth
manifold. It's known that $TM$ is trivial. If we take
$\nabla_0^{TM}$ to be the trivial connection $d^{TM}$ associated to
some global basis of $\Gamma(TM)$ and $\nabla_1^{TM}=d^{TM}+A$,
(2.5) gives us the well known Chern-Simons form  (cf. [10], [31])\be
CS(A)\triangleq\mathrm{tr}\left[A\wedge d^{TM}A+{2\over3}A\wedge A
\wedge A\right],\ee which plays a very important role in quantum
field theory and low dimensional topology [27].

In this paper, we are going to use the following theorem, which we
obtain by modifying Theorem 2.1.
\begin{theorem} Assume $a_0\neq 0$ and fix a sign when expanding $\mathrm{det}^{1\over2}(f(R_t^E)),\, 0\leq t \leq 1$. The following identity holds,
\be
\mathrm{det}^{1\over2}(f(\widetilde{R}^E))-\mathrm{det}^{1\over2}(f(R^E))=d
\int_0^1 {1\over 2}\mathrm{det}^{1\over2}(f(R_t^E))\mathrm{tr}\left[
\frac{d\nabla^E_t}{dt}\frac{f^{'}(R^E_t)}{f(R^E_t)}\right]dt.\ee
\end{theorem}
\begin{proof} Observe that $\mathrm{det}^{1\over2}(f(R^E_t))=e^{{1\over
2}\mathrm{tr}[\mathrm{ln}f(R^E_t)]}$. Hence by (2.4),

\be  \begin{split} \frac{d}{dt} \mathrm{det}^{1 \over
2}(f(R^E_t))&=\frac{d}{dt}e^{{1\over
2}\mathrm{tr}[\mathrm{ln}f(R^E_t)] } \\
&={1\over2}e^{{1\over2}\mathrm{tr}[\mathrm{ln}f(R^E_t)]}\frac{d}{dt}\mathrm{tr}[\mathrm{ln}f(R^E_t)]\\
&={1\over2}e^{{1\over2}\mathrm{tr}[\mathrm{ln}f(R^E_t)]}
d\,\mathrm{tr}\left[\frac{d\nabla^E_t}{dt}\frac{f^{'}(R^E_t)}{f(R^E_t)}\right]\\
&=d\,\left\{{1\over2}e^{{1\over2}\mathrm{tr}[\mathrm{ln}f(R^E_t)]}
\mathrm{tr}\left[\frac{d\nabla^E_t}{dt}\frac{f^{'}(R^E_t)}{f(R^E_t)}\right]\right\}\\
&=d\,\left\{{1\over2}\mathrm{det}^{1 \over 2}(f(R^E_t))
\mathrm{tr}\left[\frac{d\nabla^E_t}{dt}\frac{f^{'}(R^E_t)}{f(R^E_t)}\right]\right\}.
\end{split}\ee
Therefore (2.7) follows.

\end{proof}

\section{The Landweber-Stong forms and the Witten forms} Let $M$ be
a Riemannian manifold. Let $\nabla^{TM}$ be the associated
Levi-Civita connection on $TM$ and $R^{TM}=(\nabla^{TM})^2$ be the
curvature of $\nabla^{TM}$. Let $\widehat{A}(TM, \nabla^{TM})$ and
$L(TM, \nabla^{TM})$ be the Hirzebruch characteristic forms defined
respectively by (cf. [31])
\begin{equation}
\begin{split}
&\widehat{A}(TM, \nabla^{TM}) ={\det}^{1/2}\left({{\sqrt{-1}\over
4\pi}R^{TM} \over \sinh\left({ \sqrt{-1}\over
4\pi}R^{TM}\right)}\right), \\ &L(TM, \nabla^{TM})
={\det}^{1/2}\left({{\sqrt{-1}\over 2\pi}R^{TM} \over \tanh\left({
\sqrt{-1}\over 2\pi}R^{TM}\right)}\right).
\end{split}
\end{equation}

Let $E$, $F$ be two Hermitian vector bundles over $M$ carrying
Hermitian connections $\nabla^E$, $\nabla^F$ respectively. Let
$R^E=(\nabla^{E})^2$ (resp. $R^F=(\nabla^{F})^2$) be the curvature
of $\nabla^E$ (resp. $\nabla^F$). If we set the formal difference
$G=E-F$, then  $G$ carries an induced Hermitian connection
$\nabla^G$ in an obvious sense. We define the associated Chern
character form as (cf.[30])
\begin{equation}{\rm ch}(G,\nabla^G)={\rm
tr}\left[\exp\left({\sqrt{-1}\over 2\pi}R^E\right)\right] -{\rm
tr}\left[\exp\left({\sqrt{-1}\over
2\pi}R^F\right)\right].\end{equation} Sometimes we also need to use
the following modified Chern character (cf. [11])
\begin{equation}\widetilde{{\rm ch}}(G,\nabla^G)={\rm
tr}\left[\exp\left({\sqrt{-1}\over \pi}R^E\right)\right] -{\rm
tr}\left[\exp\left({\sqrt{-1}\over
\pi}R^F\right)\right].\end{equation}

For any complex number $t$, let
$$\Lambda_t(E)={\mathbb C}|_M+tE+t^2\Lambda^2(E)+\cdots ,
\\\ S_t(E)={\mathbb C}|_M+tE+t^2S^2(E)+\cdots$$  denote respectively
the total exterior and symmetric powers  of $E$, which live in
$K(M)[[t]].$ The following relations between these two operations
[3, Chap. 3] hold, \be S_t(E)=\frac{1}{\Lambda_{-t}(E)},\ \ \ \
 \Lambda_t(E-F)=\frac{\Lambda_t(E)}{\Lambda_t(F)}.\ee
Moreover, if $\{\omega_i \}$, $\{{\omega_j}' \}$ are formal Chern
roots for Hermitian vector bundles $E$, $F$ respectively, then
[15, Chap. 1] \be
\mathrm{ch}(\Lambda_t{(E)})=\prod\limits_i(1+e^{\omega_i}t).\ee
Therefore,  we have the following formulas for Chern character
forms, \be{\rm ch}(S_t(E) )=\frac{1}{{\rm ch}(\Lambda_{-t}(E)
)}=\frac{1}{\prod\limits_i (1-e^{\omega_i}t)}\ ,\ee \be{\rm
ch}(\Lambda_t(E-F) )=\frac{{\rm ch}(\Lambda_t(E) )}{{\rm
ch}(\Lambda_t(F)
)}=\frac{\prod\limits_i(1+e^{\omega_i}t)}{\prod\limits_j(1+e^{{\omega_j}'}t)}\
.\ee

If $W$ is a  real Euclidean vector bundle over $M$ carrying a
Euclidean connection $\nabla^W$, then its complexification
$W_\mathbb{C}=W\otimes \mathbb{C}$ is a complex vector bundle over
$M$ carrying a canonically induced Hermitian metric from that of
$W$, as well as a Hermitian connection $\nabla^{W_\mathbb{C}}$
induced from $\nabla^W$. If $E$ is a vector bundle (complex or real)
over $M$, set $\widetilde{E}=E-{{\rm dim}E}$ in $K(M)$ or $KO(M)$.

Let $q=e^{2\pi \sqrt{-1}\tau}$ with $\tau \in \mathbb{H}$, the upper
half complex plane. Set (cf. [28], [20]) \be \Theta_1(T_{\mathbb
C}M)=\bigotimes_{n=1}^\infty S_{q^n}(\widetilde{T_{\mathbb C}M})
\otimes \bigotimes_{m=1}^\infty \Lambda_{q^m}(\widetilde{T_\CC
M}),\ee

\be \Theta_2(T_{\mathbb C}M)=\bigotimes_{n=1}^\infty
S_{q^n}(\widetilde{T_{\mathbb C}M}) \otimes
\bigotimes_{m=1}^\infty \Lambda_{-q^{m-{1\over
2}}}(\widetilde{T_{\mathbb C}M}),\ee

\be \Theta_3(T_{\mathbb C}M)=\bigotimes_{n=1}^\infty
S_{q^n}(\widetilde{T_{\mathbb C}M}) \otimes
\bigotimes_{m=1}^\infty \Lambda_{q^{m-{1\over
2}}}(\widetilde{T_{\mathbb C}M}),\ee

\be \Theta(T_{\mathbb C}M)=\bigotimes_{n=1}^\infty
S_{q^n}(\widetilde{T_{\mathbb C}M}).\ee

$\Theta_1(T_{\mathbb C}M)$, $\Theta_2(T_{\mathbb C}M)$,
$\Theta_3(T_{\mathbb C}M)$ and $\Theta(T_{\mathbb C}M)$ admit
formal Fourier expansion in $q^{1/2}$ as

\be\Theta_1(T_{\mathbb C}M)=A_0(T_{\mathbb C}M)+ A_1(T_{\mathbb
C}M)q^{1/2}+\cdots,\ee

\be \Theta_2(T_{\mathbb C}M)=B_0(T_{\mathbb C}M) +B_1(T_{\mathbb
C}M)q^{1/2}+\cdots,\ee

\be \Theta_3(T_{\mathbb C}M)=C_0(T_{\mathbb C}M) +C_1(T_{\mathbb
C}M)q^{1/2}+\cdots,\ee

\be \Theta(T_{\mathbb C}M)=D_0(T_{\mathbb C}M) +D_1(T_{\mathbb
C}M)q^{1/2}+\cdots,\ee where the $A_j$'s, $B_j$'s, $C_j$'s and
$D_j$'s are elements in the semi-group formally generated by complex
vector bundles over $M$. Moreover, they carry canonically induced
connections denoted by $\nabla^{A_j}$, $\nabla^{B_j}$,
$\nabla^{C_j}$ and $\nabla^{D_j}$ respectively, and let
$\nabla^{\Theta_i(T_{\mathbb C}M)}$, $\nabla^{\Theta(T_{\mathbb
C}M)}$ be the induced connections with $q^{1/2}$-coefficients on
$\Theta_i$, $\Theta$ from the $\nabla^{A_j}$, $\nabla^{B_j}$,
$\nabla^{C_j}$ and $\nabla^{D_j}$.

If $\omega$ is a differential form on $M$, we denote by
$\omega^{(i)}$ its degree $i$ component.

\begin{definition} \be \Phi_L(\nabla^{TM}, \tau)\triangleq L(TM,\nabla^{TM})
\widetilde{\mathrm{ch}}\left(\Theta_1(T_\CC
M),\nabla^{\Theta_1(T_\CC M)}\right)\ee is called the {\bf
Landweber-Stong form} of $M$ with respect to $\nabla^{TM}$; \be
\Phi_W(\nabla^{TM}, \tau)\triangleq \widehat{A}(TM,\nabla^{TM})
\mathrm{ch}\left(\Theta_2(T_\CC M),\nabla^{\Theta_2(T_\CC
M)}\right),\ee \be \Phi_W^{'}(\nabla^{TM}, \tau)\triangleq
\widehat{A}(TM,\nabla^{TM}) \mathrm{ch}\left(\Theta_3(T_\CC
M),\nabla^{\Theta_3(T_\CC M)}\right),\ee

\be \Psi_W(\nabla^{TM}, \tau)\triangleq\widehat{A}(TM,\nabla^{TM})
\mathrm{ch}\left(\Theta(T_\CC M),\nabla^{\Theta(T_\CC
M)}\right)\ee are called the {\bf Witten forms} of $M$ with
respect to $\nabla^{TM}$.

\end{definition}

The four Jacobi theta functions are defined as follows (cf. [7]):
\be\theta(v,\tau)=2q^{1/8}\sin(\pi v)
\prod_{j=1}^\infty\left[(1-q^j)(1-e^{2\pi \sqrt{-1}v}q^j)(1-e^{-2\pi
\sqrt{-1}v}q^j)\right]\ ,\ee \be \theta_1(v,\tau)=2q^{1/8}\cos(\pi
v)
 \prod_{j=1}^\infty\left[(1-q^j)(1+e^{2\pi \sqrt{-1}v}q^j)
 (1+e^{-2\pi \sqrt{-1}v}q^j)\right]\ ,\ee
\be \theta_2(v,\tau)=\prod_{j=1}^\infty\left[(1-q^j)
 (1-e^{2\pi \sqrt{-1}v}q^{j-1/2})(1-e^{-2\pi \sqrt{-1}v}q^{j-1/2})\right]\
 ,\ee
\be \theta_3(v,\tau)=\prod_{j=1}^\infty\left[(1-q^j) (1+e^{2\pi
\sqrt{-1}v}q^{j-1/2})(1+e^{-2\pi \sqrt{-1}v}q^{j-1/2})\right]\ .\ee
They are all holomorphic functions for $(v,\tau)\in \mathbb{C \times
H}$, where $\mathbb{C}$ is the complex plane and $\mathbb{H}$ is the
upper half plane.

Let $\theta^{'}(0,\tau)=\frac{\partial}{\partial
v}\theta(v,\tau)|_{v=0}$. The {\bf Jacobi identity} [7],
$$\theta^{'}(0,\tau)=\pi \theta_1(0,\tau)
\theta_2(0,\tau)\theta_3(0,\tau)$$ holds.

Applying the Chern-Weil theory, we can express the Landweber-Stong
forms and the Witten forms in terms of theta functions and
curvatures, which looks new in the literature (cf. [20], [30]).

\begin{proposition} The following identities hold,
\be \Phi_L(\nabla^{TM}, \tau)=\mathrm{det}^{1\over
2}\left(\frac{R^{TM}}{2{\pi}^2}\frac{\theta'(0,\tau)}{\theta(\frac{R^{TM}}{2{\pi}^2},\tau)}
\frac{\theta_{1}(\frac{R^{TM}}{2{\pi}^2},\tau)}{\theta_{1}(0,\tau)}\right),
\ee

\be \Phi_W(\nabla^{TM}, \tau)=\mathrm{det}^{1\over
2}\left(\frac{R^{TM}}{4{\pi}^2}\frac{\theta'(0,\tau)}{\theta(\frac{R^{TM}}{4{\pi}^2},\tau)}
\frac{\theta_{2}(\frac{R^{TM}}{4{\pi}^2},\tau)}{\theta_{2}(0,\tau)}\right),
\ee

\be \Phi_W^{'}(\nabla^{TM}, \tau)=\mathrm{det}^{1\over
2}\left(\frac{R^{TM}}{4{\pi}^2}\frac{\theta'(0,\tau)}{\theta(\frac{R^{TM}}{4{\pi}^2},\tau)}
\frac{\theta_{3}(\frac{R^{TM}}{4{\pi}^2},\tau)}{\theta_{3}(0,\tau)}\right),
\ee

\be \Psi_W(\nabla^{TM}, \tau)=\mathrm{det}^{1\over
2}\left(\frac{R^{TM}}{4{\pi}^2}\frac{\theta'(0,\tau)}{\theta(\frac{R^{TM}}{4{\pi}^2},\tau)}
\right). \ee
\end{proposition}

Let $$ SL_2(\mathbb{Z}):= \left\{\left.\left(\begin{array}{cc}
                                      a&b\\
                                      c&d
                                     \end{array}\right)\right|a,b,c,d\in\mathbb{Z},\ ad-bc=1
                                     \right\}
                                     $$
 as usual be the modular group. Let
$$S=\left(\begin{array}{cc}
      0&-1\\
      1&0
\end{array}\right), \ \ \  T=\left(\begin{array}{cc}
      1&1\\
      0&1
\end{array}\right)$$
be the two generators of $ SL_2(\mathbb{Z})$. Their actions on
$\mathbb{H}$ are given by
$$ S:\tau\rightarrow-\frac{1}{\tau}, \ \ \ T:\tau\rightarrow\tau+1.$$

Let
$$ \Gamma_0(2)=\left\{\left.\left(\begin{array}{cc}
a&b\\
c&d
\end{array}\right)\in SL_2(\mathbb{Z})\right|c\equiv0\ \ (\rm mod \ \ 2)\right\},$$

$$ \Gamma^0(2)=\left\{\left.\left(\begin{array}{cc}
a&b\\
c&d
\end{array}\right)\in SL_2(\mathbb{Z})\right|b\equiv0\ \ (\rm mod \ \ 2)\right\}$$

$$ \Gamma_\theta=\left\{\left.\left(\begin{array}{cc}
a&b\\
c&d
\end{array}\right)\in SL_2(\mathbb{Z})\right|\left(\begin{array}{cc}
a&b\\
c&d
\end{array}\right)\equiv\left(\begin{array}{cc}
1&0\\
0&1
\end{array}\right) \mathrm{or} \left(\begin{array}{cc}
0&1\\
1&0
\end{array}\right)\ \ (\rm mod \ \ 2)\right\}$$
be the three modular subgroups of $SL_2(\mathbb{Z})$. It is known
that the generators of $\Gamma_0(2)$ are $T,ST^2ST$, the generators
of $\Gamma^0(2)$ are $STS,T^2STS$  and the generators of
$\Gamma_\theta$ are $S$, $T^2$. (cf. [7]).

If we act theta-functions by $S$ and $T$, the theta functions obey
the following transformation laws (cf. [7]), \be
\theta(v,\tau+1)=e^{\pi \sqrt{-1}\over 4}\theta(v,\tau),\ \ \
\theta\left(v,-{1}/{\tau}\right)={1\over\sqrt{-1}}\left({\tau\over
\sqrt{-1}}\right)^{1/2} e^{\pi\sqrt{-1}\tau v^2}\theta\left(\tau
v,\tau\right)\ ;\ee \be \theta_1(v,\tau+1)=e^{\pi \sqrt{-1}\over
4}\theta_1(v,\tau),\ \ \
\theta_1\left(v,-{1}/{\tau}\right)=\left({\tau\over
\sqrt{-1}}\right)^{1/2} e^{\pi\sqrt{-1}\tau v^2}\theta_2(\tau
v,\tau)\ ;\ee \be\theta_2(v,\tau+1)=\theta_3(v,\tau),\ \ \
\theta_2\left(v,-{1}/{\tau}\right)=\left({\tau\over
\sqrt{-1}}\right)^{1/2} e^{\pi\sqrt{-1}\tau v^2}\theta_1(\tau
v,\tau)\ ;\ee \be\theta_3(v,\tau+1)=\theta_2(v,\tau),\ \ \
\theta_3\left(v,-{1}/{\tau}\right)=\left({\tau\over
\sqrt{-1}}\right)^{1/2} e^{\pi\sqrt{-1}\tau v^2}\theta_3(\tau
v,\tau)\ .\ee

\begin{definition} Let $\Gamma$ be a subgroup of $SL_2(\mathbb{Z}).$ A modular form over $\Gamma$ is a holomorphic function $f(\tau)$ on $\mathbb{H}\cup
\{\infty\}$ such that for any
 $$ g=\left(\begin{array}{cc}
             a&b\\
             c&d
             \end{array}\right)\in\Gamma\ ,$$
 the following property holds
 $$f(g\tau):=f(\frac{a\tau+b}{c\tau+d})=\chi(g)(c\tau+d)^kf(\tau), $$
 where $\chi:\Gamma\rightarrow\mathbf{C}^*$ is a character of
 $\Gamma$ and $k$ is called the weight of $f$.
 \end{definition}

If $\Gamma$ is a modular subgroup, let
$\mathcal{M}_\mathbb{R}(\Gamma)$ denote the ring of modular forms
over $\Gamma$ with real Fourier coefficients. Writing simply
$\theta_j=\theta_j(0,\tau),\ 1\leq j \leq 3,$ we introduce six
explicit modular forms (cf.\ [18], [20]),
$$ \delta_1(\tau)=\frac{1}{8}(\theta_2^4+\theta_3^4), \ \ \ \
\varepsilon_1(\tau)=\frac{1}{16}\theta_2^4 \theta_3^4\ ,$$
$$\delta_2(\tau)=-\frac{1}{8}(\theta_1^4+\theta_3^4), \ \ \ \
\varepsilon_2(\tau)=\frac{1}{16}\theta_1^4 \theta_3^4\ ,$$
$$\delta_3(\tau)=\frac{1}{8}(\theta_1^4-\theta_2^4), \ \ \ \
\varepsilon_3(\tau)=-\frac{1}{16}\theta_1^4 \theta_2^4\ .$$ They
have the following Fourier expansions in $q^{1/2}$:
$$\delta_1(\tau)={1\over 4}+6q+6q^2+\cdots,\ \ \ \ \varepsilon_1(\tau)={1\over
16}-q+7q^2+\cdots\ , $$
$$\delta_2(\tau)=-{1\over 8}-3q^{1/2}-3q+\cdots,\ \ \ \
\varepsilon_2(\tau)=q^{1/2}+8q+\cdots\ ,$$

$$\delta_3(\tau)=-{1\over 8}+3q^{1/2}-3q+\cdots,\ \ \ \
\varepsilon_3(\tau)=-q^{1/2}+8q+\cdots.$$ where the
\textquotedblleft $\cdots$" terms are the higher degree terms, all
of which have integral coefficients. They also satisfy the
transformation laws (cf [18], [20], [21]), \be
\delta_2\left(-\frac{1}{\tau}\right)=\tau^2\delta_1(\tau)\ \ \ \ \ ,
\ \ \ \ \
\varepsilon_2\left(-\frac{1}{\tau}\right)=\tau^4\varepsilon_1(\tau).\ee

\be \delta_2\left(\tau+1\right)=\delta_3(\tau)\ \ \ \ \ , \ \ \ \
\ \varepsilon_2\left(\tau+1\right)=\varepsilon_3(\tau).\ee

\begin{lemma}$\mathrm{(cf. [20])}$ One has that $\delta_1(\tau)\ (resp.\ \varepsilon_1(\tau) ) $
is a modular form of weight $2 \ (resp.\ 4)$ over $\Gamma_0(2)$,
$\delta_2(\tau) \ (resp.\ \varepsilon_2(\tau))$ is a modular form
of weight $2\ (resp.\ 4)$ over $\Gamma^0(2)$, while
$\delta_3(\tau) \ (resp.\ \varepsilon_3(\tau))$ is a modular form
of weight $2\ (resp.\ 4)$ over $\Gamma_\theta(2)$ and moreover
$\mathcal{M}_\mathbb{R}(\Gamma^0(2))=\mathbb{R}[\delta_2(\tau),
\varepsilon_2(\tau)]$.
\end{lemma}

Acting the transformations $S, T$ on the Landweber-Stong forms and
the Witten forms, we have (cf. [20][21][14])
\begin{proposition} For any integer $i\geq 0$, one has that
\newline 1)\ $\left\{\Phi_L(\nabla^{TM}, \tau)\right\}^{(4i)}$ is a modular form of weight 2i over $\Gamma_0(2)$;
\newline $\left\{\Phi_W(\nabla^{TM}, \tau)\right\}^{(4i)}$ is a modular form of weight 2i over $\Gamma^0(2)$;
\newline $\left\{\Phi_W^{'}(\nabla^{TM}, \tau)\right\}^{(4i)}$ is a modular form of weight 2i over
$\Gamma_\theta$;
\newline moreover, if the first Pontrjagin form $p_1(M,\nabla^{TM})=0$, then \newline
$\left\{\Psi_W(\nabla^{TM}, \tau)\right\}^{(4i)}$ is a modular
form of weight 2i over $SL(2, \mathbb{Z})$.
\newline 2)\ the following equalities hold,
\be \left\{\Phi_L(\nabla^{TM},
 -{1}/{\tau})\right\}^{(4i)}=(2\tau)^{2i}\left\{\Phi_W(\nabla^{TM},
\tau)\right\}^{(4i)},\ee \be \Phi_W(\nabla^{TM},
\tau+1)=\Phi_W'(\nabla^{TM}, \tau).\ee
\end{proposition}

The modularities in Proposition 3.2 have some interesting
applications. For example, let $M$ be 12-dimensional and $i=3$.
$\left\{\Phi_L(\nabla^{TM}, \tau)\right\}^{(12)}$ is a modular form
of weight 6 over $\Gamma_0(2)$, $\left\{\Phi_W(\nabla^{TM},
\tau)\right\}^{(12)}$ is a modular form of weight 6 over
$\Gamma^0(2)$ and $$ \left\{\Phi_L(\nabla^{TM},
-{1}/{\tau})\right\}^{(12)}=(2\tau)^{6}\left\{\Phi_W(\nabla^{TM},
 \tau)\right\}^{(12)}.$$ Then by Lemma 3.1, we have \be
\left\{\Phi_W(\nabla^{TM},
\tau)\right\}^{(12)}=h_0(8\delta_2)^3+h_1
(8\delta_2)\epsilon_2,\ee and by (3.32) and (3.34), \be
\left\{\Phi_L(\nabla^{TM},
\tau)\right\}^{(12)}=2^6[h_0(8\delta_1)^3+h_1
(8\delta_1)\epsilon_1], \ee where by comparing the
$q^{1/2}$-expansion coefficients in (3.36),
$h_0=-\{\widehat{A}(TM, \nabla^{TM})\}^{(12)}$ and
$h_1=\{60\widehat{A}(TM, \nabla^{TM})+\widehat{A}(TM,
\nabla^{TM})\mathrm{ch}(T_\CC M, \nabla^{T_\CC M}\}^{(12)}$. Then
comparing the constant coefficients of the $q$-expansions of both
sides of (3.37), one obtains that $\{L(TM,
\nabla^{TM})\}^{(12)}=2^3(2^6h_0+h_1)$, consequently \be \{L(TM,
\nabla^{TM})\}^{(12)}=\{8\widehat{A}(TM,
\nabla^{TM})\mathrm{ch}(T_\CC M, \nabla^{T_\CC
M})-32\widehat{A}(TM, \nabla^{TM})\}^{(12)},\ee which is just the
gravitational anomaly cancellation formula derived by
Alvarez-Gaum\'e and Witten in [1] from very nontrivial
computations. Liu [20] generalizes the miraculous cancellation
formula (3.38) to arbitrary $8k+4$ dimensional smooth manifolds by
developing modular invariance properties of characteristic forms.
Formulas of this type have interesting applications in the study
of divisibility and congruence phenomena for characteristic
numbers. We refer interested readers to [12][13][20].

Let $M$ be a $4k$ dimensional closed oriented smooth manifold and
$[M]$ is the fundamental class. $\phi_L(M,
\tau)\triangleq\langle\Phi_L(\nabla^{TM}, \tau), [M]\rangle$ is a
modular form of weight $2k$ over $\Gamma_0(2)$ with integral
$q$-expansion coefficients and called the {\bf Landweber-Stong
genus} of $M$. $\phi_W(M, \tau)\triangleq\langle\Phi_W(\nabla^{TM},
\tau), [M]\rangle$ and $\phi_W^{'}(M,
\tau)\triangleq\langle\Phi_W^{'}(\nabla^{TM}, \tau), [M]\rangle$ are
modular forms of weight $2k$ over $\Gamma^0(2)$ and $\Gamma_\theta$
respectively. If the first Pontrajagin class of $M$ vanishes, then
$\psi_W(M, \tau)\triangleq\langle \Psi_W(\nabla^{TM}, \tau),
[M]\rangle$ is a modular form over $SL(2, \mathbb{Z})$. They are
modular forms with rational $q$-expansion coefficients and called
the {\bf Witten genera} of $M$. These are all examples of the {\bf
elliptic genera}, which were first defined by Ochanine [24].

Let $d_s$ be the signature operator of $M$. According to the
Atiyah-Singer index theorem, we can express the Landweber-Stong
genus analytically by the index of the twisted signature operator
as
$$\phi_L(M, \tau)={\mathrm {Ind}}(d_s \otimes \Theta_1(T_{\mathbb C}M)).
$$ Moreover let $M$ be spin and $D$ be the Atiyah-Singer Dirac
operator over $M$. The Witten genera can also be analytically
expressed by the indices of the twisted Dirac operators as
followings,
$$\phi_W(M, \tau)={\mathrm {Ind}}(D\otimes \Theta_2(T_{\mathbb C}M)), $$
$$\phi_W^{'}(M, \tau)={\mathrm {Ind}}(D\otimes \Theta_3(T_{\mathbb
C}M)),$$
$$ \psi_W(M, \tau)={\mathrm {Ind}}(D\otimes
\Theta(T_{\mathbb C}M)).$$ Heuristically, these twisted operators
are viewed as elliptic operators on the smooth loop space $LM$ from
the view point of string theory. See details in Section 5. The
elliptic operators $d_s\otimes \Theta_1(T_{\mathbb C}M), D\otimes
\Theta_2(T_{\mathbb C}M)$ and $D\otimes \Theta_3(T_{\mathbb C}M)$
are all rigid according to the Witten rigidity theorem, which is
proved by Bott-Taubes and Liu ([6], [22]).

\section{Transgressed forms and modularities}

Consider the following functions defined on $\CC \times
\mathbb{H}$,
$$f_{\Phi_L}(z, \tau)=2z\frac{\theta'(0,\tau)}{\theta(2z,\tau)}
\frac{\theta_{1}(2z,\tau)}{\theta_{1}(0,\tau)},$$
$$f_{\Phi_W}(z,
\tau)=z\frac{\theta'(0,\tau)}{\theta(z,\tau)}
\frac{\theta_{2}(z,\tau)}{\theta_{2}(0,\tau)},$$
$$f_{\Phi_W^{'}}(z,
\tau)=z\frac{\theta'(0,\tau)}{\theta(z,\tau)}
\frac{\theta_{3}(z,\tau)}{\theta_{3}(0,\tau)},$$
$$f_{\Psi_W}(z,
\tau)=z\frac{\theta'(0,\tau)}{\theta(z,\tau)}.$$ For the $4k$
dimensional manifold $M$, by Proposition 3.1, we have
$$\Phi_L(\nabla^{TM}, \tau)=\mathrm{det}^{1\over 2}\left(f_{\Phi_L}\left(\frac{R^{TM}}{4\pi^2}, \tau\right)\right),$$
$$\Phi_W(\nabla^{TM}, \tau)=\mathrm{det}^{1\over 2}\left(f_{\Phi_W}\left(\frac{R^{TM}}{4\pi^2}, \tau\right)\right),$$
$$\Phi_W^{'}(\nabla^{TM}, \tau)=\mathrm{det}^{1\over 2}\left(f_{\Phi_W^{'}}\left(\frac{R^{TM}}{4\pi^2}, \tau\right)\right),$$
$$\Psi_W(\nabla^{TM}, \tau)=\mathrm{det}^{1\over 2}\left(f_{\Psi_W}\left(\frac{R^{TM}}{4\pi^2}, \tau\right)\right).$$

From now on, let $M$ be a $4k-1$ dimensional smooth manifold. Let
$\nabla^{TM}_i, i=0,1$ be two connections on $TM$ and $R^{TM}_i,
i=0,1$ be their curvatures respectively. Let
$\nabla_t^{TM}=(1-t)\nabla_0^{TM}+t\nabla_1^{TM}$ and $R_t^{TM}$ be
the corresponding curvature. Let $A=\nabla_1-\nabla_0 \in
\Omega^1(M, \mathrm{End}(TM)).$

By Theorem 2.2, one has \be
\begin{split} &\mathrm{det}^{1\over
2}\left(f_{\Phi_L}\left(\frac{R_1^{TM}}{4\pi^2},
\tau\right)\right)-\mathrm{det}^{1\over
2}\left(f_{\Phi_L}\left(\frac{R_0^{TM}}{4\pi^2},
\tau\right)\right)\\
=&d \int_0^1 {\frac{1} {8\pi^2}}\mathrm{det}^{\frac{1}
{2}}\left(f_{\Phi_L}\left(\frac{R_t^{TM}}{4\pi^2},
\tau\right)\right)\mathrm{tr}\left[A\frac{f_{\Phi_L}'(\frac{R_t^{TM}}{4\pi^2}, \tau)}{f_{\Phi_L}(\frac{R_t^{TM}}{4\pi^2}, \tau)}\right]dt.\\
\end{split}\ee
We define \be \begin{split} &\mathrm{CS}\Phi_L(\nabla_0^{TM},
\nabla_1^{TM},
\tau)\\
&\triangleq {\frac{1} {4\pi^2}}\int_0^1 \Phi_L(\nabla_t^{TM},
\tau)\mathrm{tr}\left[A\left(\frac{1}{\frac{R_t^{TM}}{2\pi^2}}-\frac{\theta'(\frac{R_t^{TM}}{2\pi^2},
\tau)}{\theta(\frac{R_t^{TM}}{2\pi^2},
\tau)}+\frac{\theta_1'(\frac{R_t^{TM}}{2\pi^2},
\tau)}{\theta_1(\frac{R_t^{TM}}{2\pi^2}, \tau)}\right)\right]dt
,\end{split}\ee which is in $\Omega^{\mathrm{odd}}(M, \CC)[[q]]$.
Since $M$ is $4k-1$ dimensional,
$\left\{\mathrm{CS}\Phi_L(\nabla_0^{TM}, \nabla_1^{TM},
\tau)\right\}^{(4k-1)}$ represents an element in $H^{4k-1}(M,
\mathbb{C})[[q]]$.

Similarly, we can compute the transgressed forms for $\Phi_W,
\Phi_W'$ and $\Psi_W$ respectively and define \be \begin{split}
&\mathrm{CS}\Phi_W(\nabla_0^{TM}, \nabla_1^{TM},
\tau)\\
&\triangleq{\frac{1} {8\pi^2}}\int_0^1 \Phi_W(\nabla_t^{TM},
\tau)\mathrm{tr}\left[A\left(\frac{1}{\frac{R_t^{TM}}{4\pi^2}}-\frac{\theta'(\frac{R_t^{TM}}{4\pi^2},
\tau)}{\theta(\frac{R_t^{TM}}{4\pi^2},
\tau)}+\frac{\theta_2'(\frac{R_t^{TM}}{4\pi^2},
\tau)}{\theta_2(\frac{R_t^{TM}}{4\pi^2}, \tau)}\right)\right]
dt,\end{split}\ee

\be \begin{split}&\mathrm{CS}\Phi_W'(\nabla_0^{TM}, \nabla_1^{TM},
\tau)\\
&\triangleq{\frac{1} {8\pi^2}}\int_0^1 \Phi_W'(\nabla_t^{TM},
\tau)\mathrm{tr}\left[A\left(\frac{1}{\frac{R_t^{TM}}{4\pi^2}}-\frac{\theta'(\frac{R_t^{TM}}{4\pi^2},
\tau)}{\theta(\frac{R_t^{TM}}{4\pi^2},
\tau)}+\frac{\theta_3'(\frac{R_t^{TM}}{4\pi^2},
\tau)}{\theta_3(\frac{R_t^{TM}}{4\pi^2}, \tau)}\right)\right]
dt,\end{split}\ee and

\be \begin{split} &\mathrm{CS}\Psi_W(\nabla_0^{TM}, \nabla_1^{TM},
\tau)\\
&\triangleq{\frac{1} {8\pi^2}}\int_0^1 \Psi_W(\nabla_t^{TM},
\tau)\mathrm{tr}\left[A\left(\frac{1}{\frac{R_t^{TM}}{4\pi^2}}-\frac{\theta'(\frac{R_t^{TM}}{4\pi^2},
\tau)}{\theta(\frac{R_t^{TM}}{4\pi^2}, \tau)}\right)\right]
dt,\end{split}\ee which lie in $\Omega^{\mathrm{odd}}(M,
\mathbb{C})[[q^{1\over 2}]]$ and their top components represent
elements in $H^{4k-1}(M, \mathbb{C})[[q^{1\over2}]]$.

\begin{remark} In (4.2) to (4.5) and in the following, to make sense,
${1\over z} -\frac{\theta'(z, \tau)}{\theta(z, \tau)}$ should be
understood as the z-expansion of $ \frac{[\theta(z,
\tau)-z\theta'(z, \tau)]/z^2}{\theta(z, \tau)/z},$ where both the
nominator and the denominator have nonzero constant terms in their
z-expansions.
\end{remark}

Equality (4.1) and the modular invariance properties of
$\mathrm{det}^{1\over
2}\left(f_{\Phi_L}\left(\frac{R_1^{TM}}{4\pi^2},
\tau\right)\right)$, $\mathrm{det}^{1\over
2}\left(f_{\Phi_L}\left(\frac{R_0^{TM}}{4\pi^2},
\tau\right)\right)$ are not enough to guarantee that
$\mathrm{CS}\Phi_L(\nabla_0^{TM}, \nabla_1^{TM}, \tau)$ is a
modular form. Actually ($\mathrm{CS}\Phi_L(\nabla_0^{TM},
\nabla_1^{TM}, \tau)$+ any\ closed\  form) will also satisfy
(4.1). This also true for other transgressed forms (4.3)-(4.5).
However we do have the following results,
\begin{theorem}Let $M$ be a $4k-1$ dimensional smooth manifold and $\nabla^{TM}_0, \nabla^{TM}_1$ be two connections on $TM$,
then for integer $i, 1\leq i \leq k$, we have
\newline 1)\ $\left\{\mathrm{CS}\Phi_L(\nabla_0^{TM}, \nabla_1^{TM},
\tau)\right\}^{(4i-1)}$ is a modular form of weight 2i over
$\Gamma_0(2)$;
\newline $\left\{\mathrm{CS}\Phi_W(\nabla_0^{TM}, \nabla_1^{TM},
\tau)\right\}^{(4i-1)}$ is a modular form of weight 2i over
$\Gamma^0(2)$;
\newline $\left\{\mathrm{CS}\Phi_W^{'}(\nabla_0^{TM}, \nabla_1^{TM},
\tau)\right\}^{(4i-1)}$ is a modular form of weight 2i over
$\Gamma_\theta$;
\newline 2)\ the following equalities hold,
$$\left\{\mathrm{CS}\Phi_L(\nabla_0^{TM}, \nabla_1^{TM},
-{1}/{\tau})\right\}^{(4i-1)}=(2\tau)^{2i}\left\{\mathrm{CS}\Phi_W(\nabla_0^{TM},
\nabla_1^{TM}, \tau)\right\}^{(4i-1)},$$
$$\mathrm{CS}\Phi_W(\nabla_0^{TM}, \nabla_1^{TM},
\tau+1)=\mathrm{CS}\Phi_W'(\nabla_0^{TM}, \nabla_1^{TM}, \tau).$$
\end{theorem}

\begin{proof} Differentiating the transformation formulas (3.28) to
(3.31), we obtain that \be \begin{split} &\theta'(v,\tau+1)=e^{\pi
\sqrt{-1}\over 4}\theta'(v,\tau), \\
&\theta'\left(v,-{1}/{\tau}\right)={1\over\sqrt{-1}}\left({\tau\over
\sqrt{-1}}\right)^{1/2} e^{\pi\sqrt{-1}\tau v^2}(2\pi\sqrt{-1}\tau
v\theta\left(\tau v,\tau\right)+\tau\theta'(\tau v, \tau));\\
&\theta_1'(v,\tau+1)=e^{\pi \sqrt{-1}\over
4}\theta_1'(v,\tau),\\
&\theta_1'\left(v,-{1}/{\tau}\right)=\left({\tau\over
\sqrt{-1}}\right)^{1/2} e^{\pi\sqrt{-1}\tau v^2}(2\pi\sqrt{-1}\tau
v\theta_2\left(\tau v,\tau\right)+\tau\theta_2'(\tau v, \tau));\\
&\theta_2'(v,\tau+1)=\theta_3'(v,\tau),\\
&\theta_2'\left(v,-{1}/{\tau}\right)=\left({\tau\over
\sqrt{-1}}\right)^{1/2} e^{\pi\sqrt{-1}\tau v^2}(2\pi\sqrt{-1}\tau
v\theta_1\left(\tau v,\tau\right)+\tau\theta_1'(\tau v, \tau));\\
&\theta_3'(v,\tau+1)=\theta_2'(v,\tau),\\
&\theta_3'\left(v,-{1}/{\tau}\right)=\left({\tau\over
\sqrt{-1}}\right)^{1/2} e^{\pi\sqrt{-1}\tau v^2}(2\pi\sqrt{-1}\tau
v\theta_3\left(\tau v,\tau\right)+\tau\theta_3'(\tau v,
\tau)).\end{split}\ee Therefore \be
\theta'\left(0,-{1}/{\tau}\right)={1\over\sqrt{-1}}\left({\tau\over
\sqrt{-1}}\right)^{1/2} \tau\theta'(0, \tau).\ee By (3.28), (3.29)
and (4.7), we have \be \begin{split}
&2z\frac{\theta'(0,-{1}/{\tau})}{\theta(2z,-{1}/{\tau})}
\frac{\theta_{1}(2z,-{1}/{\tau})}{\theta_{1}(0,-{1}/{\tau})}\\
=&2z\frac{{1\over\sqrt{-1}}\left({\tau\over \sqrt{-1}}\right)^{1/2}
\tau\theta'(0, \tau)}{{1\over\sqrt{-1}}\left({\tau\over
\sqrt{-1}}\right)^{1/2} e^{\pi\sqrt{-1}\tau (2z)^2}\theta\left(2\tau
z,\tau\right)} \frac{\left({\tau\over \sqrt{-1}}\right)^{1/2}
e^{\pi\sqrt{-1}\tau (2z)^2}\theta_2(2\tau z,\tau)}{\left({\tau\over
\sqrt{-1}}\right)^{1/2}\theta_2(0,\tau)}\\
=&2\tau z\frac{\theta'(0, \tau)}{\theta\left(2\tau z,\tau\right)}
\frac{\theta_2(2\tau z,\tau)}{\theta_2(0,\tau)}.\\
\end{split} \ee
By (3.28), (3.29) and (4.6), one has \be \begin{split}
&\frac{1}{2z}-\frac{\theta'(2z, -{1}/{\tau})}{\theta(2z,
-{1}/{\tau})}+\frac{\theta_1'(2z, -{1}/{\tau})}{\theta_1(2z,
-{1}/{\tau})}\\
=&\frac{1}{2z}-\frac{{1\over\sqrt{-1}}\left({\tau\over
\sqrt{-1}}\right)^{1/2} e^{\pi\sqrt{-1}\tau
(2z)^2}(2\pi\sqrt{-1}(2\tau z)\theta\left(2\tau
z,\tau\right)+\tau\theta'(2\tau z,
\tau))}{{1\over\sqrt{-1}}\left({\tau\over \sqrt{-1}}\right)^{1/2}
e^{\pi\sqrt{-1}\tau (2z)^2}\theta\left(2\tau z,\tau\right)}\\
&+\frac{\left({\tau\over \sqrt{-1}}\right)^{1/2} e^{\pi\sqrt{-1}\tau
(2z)^2}(2\pi\sqrt{-1}(2\tau z)\theta_2\left(2\tau
z,\tau\right)+\tau\theta_2'(2\tau z, \tau))}{\left({\tau\over
\sqrt{-1}}\right)^{1/2} e^{\pi\sqrt{-1}\tau (2z)^2}\theta_2(2\tau
z,\tau)}\\
=&\frac{1}{2z}-2\pi\sqrt{-1}(2\tau z)-\tau \frac{\theta'\left(2\tau
z,\tau\right)}{\theta\left(2\tau z,\tau\right)}+2\pi\sqrt{-1}(2\tau
z)+\tau \frac{\theta_2'\left(2\tau
z,\tau\right)}{\theta_2\left(2\tau
z,\tau\right)}\\
=&\tau\left(\frac{1}{2\tau z}-\frac{\theta'\left(2\tau
z,\tau\right)}{\theta\left(2\tau z,\tau\right)}+
\frac{\theta_2'\left(2\tau z,\tau\right)}{\theta_2\left(2\tau
z,\tau\right)} \right).
\end{split} \ee
Therefore \be \begin{split} &\mathrm{CS}\Phi_L(\nabla_0^{TM},
\nabla_1^{TM},
-{1}/{\tau})\\
=&{\frac{1} {4\pi^2}}\int_0^1 \mathrm{det}^{1\over
2}\left(f_{\Phi_L}\left(\frac{R_t^{TM}}{4\pi^2},
-{1}/{\tau}\right)\right)\mathrm{tr}\left[A\left(\frac{1}{\frac{R_t^{TM}}{2\pi^2}}-\frac{\theta'(\frac{R_t^{TM}}{2\pi^2},
-{1}/{\tau})}{\theta(\frac{R_t^{TM}}{2\pi^2},
-{1}/{\tau})}+\frac{\theta_1'(\frac{R_t^{TM}}{2\pi^2},
-{1}/{\tau})}{\theta_1(\frac{R_t^{TM}}{2\pi^2},
-{1}/{\tau})}\right)\right]dt\\
=& {\frac{2\tau} {8\pi^2}}\int_0^1\mathrm{det}^{1\over
2}\left(f_{\Phi_W}\left(\frac{2\tau R_t^{TM}}{4\pi^2},
{\tau}\right)\right)\mathrm{tr}\left[A\left(\frac{1}{\frac{2\tau
R_t^{TM}}{4\pi^2}}-\frac{\theta'(\frac{2\tau R_t^{TM}}{4\pi^2},
{\tau})}{\theta(\frac{2\tau R_t^{TM}}{4\pi^2},
{\tau})}+\frac{\theta_2'(\frac{2\tau R_t^{TM}}{4\pi^2},
{\tau})}{\theta_2(\frac{2\tau R_t^{TM}}{4\pi^2},
{\tau})}\right)\right]dt. \end{split}\ee Note that the $(4i-1)$
component of the right hand side of (4.10) consists of terms like
$$\int_0^1 \mathrm{tr}\left[\left(\frac{2\tau
R_t^{TM}}{4\pi^2}\right)^{p_1}\right]\cdots\mathrm{tr}\left[\left(\frac{2\tau
R_t^{TM}}{4\pi^2}\right)^{p_s}\right]\mathrm{tr}\left[A\left(\frac{2\tau
R_t^{TM}}{4\pi^2}\right)^q \right]dt,$$ where
$2p_1+\cdots+2p_s+2q+1=4i-1$, i.e. $p_1+\cdots+p_s+q=2i-1$ due to
the fact that $A\in \Omega^1(M, \mathrm{End}(TM))$ and
$R_t^{TM}\in \Omega^2(M, \mathrm{End}(TM))$. Hence we have \be
\begin{split} &\left\{\mathrm{CS}\Phi_L(\nabla_0^{TM},
\nabla_1^{TM},
-{1}/{\tau})\right\}^{(4i-1)}\\
=&\left\{{\frac{2\tau} {8\pi^2}}\int_0^1 \mathrm{det}^{1\over
2}\left(f_{\Phi_W}\left(\frac{2\tau R_t^{TM}}{4\pi^2},
{\tau}\right)\right)\mathrm{tr}\left[A\left(\frac{1}{\frac{2\tau
R_t^{TM}}{4\pi^2}}-\frac{\theta'(\frac{2\tau R_t^{TM}}{4\pi^2},
{\tau})}{\theta(\frac{2\tau R_t^{TM}}{4\pi^2},
{\tau})}+\frac{\theta_2'(\frac{2\tau R_t^{TM}}{4\pi^2},
{\tau})}{\theta_2(\frac{2\tau R_t^{TM}}{4\pi^2},
{\tau})}\right)\right]dt\right\}^{(4i-1)}\\
=&(2\tau)^{1+2i-1}\left\{{\frac{1} {8\pi^2}}\int_0^1
\Phi_W(\nabla_t^{TM}, \tau)\mathrm{tr}\left[A\left(\frac{1}{\frac{
R_t^{TM}}{4\pi^2}}-\frac{\theta'(\frac{ R_t^{TM}}{4\pi^2},
{\tau})}{\theta(\frac{R_t^{TM}}{4\pi^2},
{\tau})}+\frac{\theta_2'(\frac{ R_t^{TM}}{4\pi^2},
{\tau})}{\theta_2(\frac{ R_t^{TM}}{4\pi^2},
{\tau})}\right)\right]dt\right\}^{(4i-1)}\\
=&(2\tau)^{2i}\left\{\mathrm{CS}\Phi_W(\nabla_0^{TM},
\nabla_1^{TM}, {\tau})\right\}^{(4i-1)}.\end{split}\ee

Similarly applying the transformation laws (3.28) to (3.31) and
(4.6), we can show that \be
\begin{split} &\mathrm{CS}\Phi_L(\nabla_0^{TM}, \nabla_1^{TM},
\tau+1)=\mathrm{CS}\Phi_L(\nabla_0^{TM}, \nabla_1^{TM},
\tau),\\
&\{\mathrm{CS}\Phi_W(\nabla_0^{TM}, \nabla_1^{TM},
-{1}/{\tau})\}^{(4i-1)}=\left(\frac{\tau}{2}\right)^{2i}\{\mathrm{CS}\Phi_L(\nabla_0^{TM},
\nabla_1^{TM}, \tau)\}^{(4i-1)},\\
&\mathrm{CS}\Phi_W(\nabla_0^{TM}, \nabla_1^{TM},
\tau+1)=\mathrm{CS}\Phi_W'(\nabla_0^{TM}, \nabla_1^{TM}, \tau),\\
&\{\mathrm{CS}\Phi_W'(\nabla_0^{TM}, \nabla_1^{TM},
-{1}/{\tau})\}^{(4i-1)}=\tau^{2i}\{\mathrm{CS}\Phi_W'(\nabla_0^{TM},
\nabla_1^{TM}, \tau)\}^{(4i-1)},\\
&\mathrm{CS}\Phi_W'(\nabla_0^{TM}, \nabla_1^{TM},
\tau+1)=\mathrm{CS}\Phi_W(\nabla_0^{TM}, \nabla_1^{TM}, \tau).\\
\end{split}\ee

Acting $ST^2ST$ to $\mathrm{CS}\Phi_L(\nabla_0^{TM},
\nabla_1^{TM}, \tau)$, we can see from (4.11) and (4.12) that \be
\begin{split} &\left\{\mathrm{CS}\Phi_L(\nabla_0^{TM},
\nabla_1^{TM},
ST^2ST\tau)\right\}^{(4i-1)}\\
=&\left\{\mathrm{CS}\Phi_L(\nabla_0^{TM}, \nabla_1^{TM},
S(T^2ST\tau))\right\}^{(4i-1)}\\
=&(2T^2ST\tau)^{2i}\left\{\mathrm{CS}\Phi_W(\nabla_0^{TM},
\nabla_1^{TM},
(T^2ST\tau))\right\}^{(4i-1)}\\
=&(2T^2ST\tau)^{2i}\left\{\mathrm{CS}\Phi_W'(\nabla_0^{TM},
\nabla_1^{TM},
(TST\tau))\right\}^{(4i-1)}\\
=&(2T^2ST\tau)^{2i}\left\{\mathrm{CS}\Phi_W(\nabla_0^{TM},
\nabla_1^{TM},
(ST\tau))\right\}^{(4i-1)}\\
=&(2T^2ST\tau)^{2i}\left(\frac{T\tau}{2}\right)^{2i}\left\{\mathrm{CS}\Phi_L(\nabla_0^{TM},
\nabla_1^{TM},
(T\tau))\right\}^{(4i-1)}\\
=&(2\frac{2\tau+1}{\tau+1})^{2i}\left(\frac{\tau+1}{2}\right)^{2i}\left\{\mathrm{CS}\Phi_L(\nabla_0^{TM},
\nabla_1^{TM},
(T\tau))\right\}^{(4i-1)}\\
=&(2\tau+1)^{2i}\left\{\mathrm{CS}\Phi_L(\nabla_0^{TM},
\nabla_1^{TM}, \tau)\right\}^{(4i-1)}. \end{split} \ee  Note that
$ST^2ST\tau=-\frac{\tau+1}{2\tau+1}$. By the first equality in
(4.12) and (4.13), we see that
$\left\{\mathrm{CS}\Phi_L(\nabla_0^{TM}, \nabla_1^{TM},
\tau)\right\}^{(4i-1)}$ is modularly invariant under the actions
of $T$ and $ST^2ST$, which form a basis for $\Gamma_0(2)$. Thus
$\left\{\mathrm{CS}\Phi_L(\nabla_0^{TM}, \nabla_1^{TM},
\tau)\right\}^{(4i-1)}$ is a modular form of weight $2i$ over
$\Gamma_0(2).$

We can similarly show that
$\left\{\mathrm{CS}\Phi_W(\nabla_0^{TM}, \nabla_1^{TM},
\tau)\right\}^{(4i-1)}$ is a modular form of weight $2i$ over
$\Gamma^0(2)$ and $\left\{\mathrm{CS}\Phi_W^{'}(\nabla_0^{TM},
\nabla_1^{TM}, \tau)\right\}^{(4i-1)}$ is a modular form of weight
$2i$ over $\Gamma_\theta$.
\end{proof}

Let's take a look at a concrete example. Let $M$ be a compact
oriented smooth 3-dimensional manifold. We have \be
\begin{split} &\mathrm{CS}\Phi_L(\nabla_0^{TM}, \nabla_1^{TM},
\tau)\\
=& {\frac{1} {4\pi^2}}\int_0^1 \Phi_L(\nabla^{TM},
\tau)\mathrm{tr}\left[A\left(\frac{1}{\frac{R_t^{TM}}{2\pi^2}}-\frac{\theta'(\frac{R_t^{TM}}{2\pi^2},
\tau)}{\theta(\frac{R_t^{TM}}{2\pi^2},
\tau)}+\frac{\theta_1'(\frac{R_t^{TM}}{2\pi^2},
\tau)}{\theta_1(\frac{R_t^{TM}}{2\pi^2},
\tau)}\right)\right]dt\\
=& {\frac{1} {4\pi^2}}\int_0^1
\mathrm{tr}\left[A\left(\frac{1}{\frac{R_t^{TM}}{2\pi^2}}-\frac{\theta'(\frac{R_t^{TM}}{2\pi^2},
\tau)}{\theta(\frac{R_t^{TM}}{2\pi^2},
\tau)}+\frac{\theta_1'(\frac{R_t^{TM}}{2\pi^2},
\tau)}{\theta_1(\frac{R_t^{TM}}{2\pi^2},
\tau)}\right)\right]dt\\
=&\frac{1} {8\pi^4}\left.\frac{\partial}{\partial
z}\left(\frac{1}{z}-\frac{\theta'(z, \tau)}{\theta(z,
\tau)}+\frac{\theta_1'(z, \tau)}{\theta_1(z,
\tau)}\right)\right|_{z=0}\int_0^1\mathrm{tr}[AR^{TM}_t]dt,\end{split}
\ee where the second equality holds because the dimension of $M$
is only 3 while $$\int_0^1
\mathrm{tr}\left[A\left(\frac{1}{\frac{R_t^{TM}}{2\pi^2}}-\frac{\theta'(\frac{R_t^{TM}}{2\pi^2},
\tau)}{\theta(\frac{R_t^{TM}}{2\pi^2},
\tau)}+\frac{\theta_1'(\frac{R_t^{TM}}{2\pi^2},
\tau)}{\theta_1(\frac{R_t^{TM}}{2\pi^2}, \tau)}\right)\right]dt$$
gives differential forms of degree greater than or equal to 3.

But \be
\begin{split}&\int_0^1\mathrm{tr}[AR^{TM}_t]dt\\
=&\int_0^1\mathrm{tr}[A((1-t)\nabla_0^{TM}+t\nabla_1^{TM})^2]dt\\
=&\int_0^1\mathrm{tr}[A((1-t)^2(\nabla_0^{TM})^2+(1-t)t[\nabla_0^{TM},\nabla_1^{TM}]+t^2(\nabla_1^{TM})^2)]dt\\
=&\mathrm{tr}\left[A({1\over
3}(\nabla_0^{TM})^2+{1\over6}[\nabla_0^{TM},\nabla_1^{TM}]+{1\over
3}(\nabla_1^{TM})^2)\right]\\
=&{1\over
3}\mathrm{tr}\left[A((\nabla_0^{TM}-\nabla_1^{TM})^2+{3\over2}[\nabla_0^{TM},\nabla_1^{TM}])\right]\\
=&{1\over
2}\mathrm{tr}\left[A[\nabla_0^{TM},\nabla_1^{TM}]+{2\over3}A\wedge
A \wedge A\right].
\end{split}\ee
Note that $\left.\frac{\partial}{\partial
z}\left(\frac{1}{z}-\frac{\theta'(z, \tau)}{\theta(z,
\tau)}+\frac{\theta_1'(z, \tau)}{\theta_1(z,
\tau)}\right)\right|_{z=0}$ is a modular form of weight 2 over
$\Gamma_0(2)$. Then by Lemma 3.1, it should be a scalar multiple
of $\delta_1(\tau)$. By direct computations, we can see that
$\left.\frac{\partial}{\partial
z}\left(\frac{1}{z}-\frac{\theta'(z, \tau)}{\theta(z,
\tau)}+\frac{\theta_1'(z, \tau)}{\theta_1(z,
\tau)}\right)\right|_{z=0}=-{2\over 3}\pi^2+O(q).$ So
$$\left.\frac{\partial}{\partial
z}\left(\frac{1}{z}-\frac{\theta'(z, \tau)}{\theta(z,
\tau)}+\frac{\theta_1'(z, \tau)}{\theta_1(z,
\tau)}\right)\right|_{z=0}=-{8\over 3}\pi^2 \delta_1(\tau).$$ Thus
we have \be \mathrm{CS}\Phi_L(\nabla_0^{TM}, \nabla_1^{TM},
\tau)=-{1\over
{6\pi^2}}\delta_1(\tau)\mathrm{tr}\left[A[\nabla_0^{TM},\nabla_1^{TM}]+{2\over3}A\wedge
A \wedge A\right].\ee

Similarly, we obtain that \be \mathrm{CS}\Phi_W(\nabla_0^{TM},
\nabla_1^{TM}, \tau)=-{1\over
{24\pi^2}}\delta_2(\tau)\mathrm{tr}\left[A[\nabla_0^{TM},\nabla_1^{TM}]+{2\over3}A\wedge
A \wedge A\right],\ee and \be \mathrm{CS}\Phi_W'(\nabla_0^{TM},
\nabla_1^{TM}, \tau)=-{1\over
{24\pi^2}}\delta_3(\tau)\mathrm{tr}\left[A[\nabla_0^{TM},\nabla_1^{TM}]+{2\over3}A\wedge
A \wedge A\right].\ee

Note that $TM$ is trivial. Let's take $\nabla_0^{TM}$ to be the
trivial connection $d^{TM}$ associated to some global basis of
$\Gamma(TM)$ and $\nabla_1^{TM}=d^{TM}+A$. We therefore have
\begin{proposition} When $M$ is a compact
oriented smooth 3-dimensional manifold, the following identities
hold, \be \mathrm{CS}\Phi_L(d^{TM}, d^{TM}+A, \tau)=-{1\over
{6\pi^2}}\delta_1(\tau)\mathrm{CS}(A),\ee \be
\mathrm{CS}\Phi_W(d^{TM}, d^{TM}+A, \tau)=-{1\over
{24\pi^2}}\delta_2(\tau)\mathrm{CS}(A),\ee \be
\mathrm{CS}\Phi_W'(d^{TM}, d^{TM}+A, \tau)=-{1\over
{24\pi^2}}\delta_3(\tau)\mathrm{CS}(A).\ee \end{proposition}

Very similar to the application after Proposition 3.2, the
modularities in Theorem 4.1 also imply some relations among
transgressed forms, which might be viewed as anomaly cancellation
formulas for odd dimensional manifolds. For example, let $M$ be 11
dimensional and $i=3$. We also similarly have that
$\left\{\mathrm{CS}\Phi_L(\nabla_0^{TM}, \nabla_1^{TM},
\tau)\right\}^{(11)}$ is a modular form of weight 6 over
$\Gamma_0(2)$, $\left\{\mathrm{CS}\Phi_W(\nabla_0^{TM},
\nabla_1^{TM}, \tau)\right\}^{(11)}$ is a modular form of weight 6
over $\Gamma^0(2)$ and $$ \left\{\mathrm{CS}\Phi_L(\nabla_0^{TM},
\nabla_1^{TM},
-{1}/{\tau})\right\}^{(11)}=(2\tau)^{6}\left\{\mathrm{CS}\Phi_W(\nabla_0^{TM},
\nabla_1^{TM}, \tau)\right\}^{(11)}.$$ Then still by Lemma 3.1, we
have \be \left\{\mathrm{CS}\Phi_W(\nabla_0^{TM}, \nabla_1^{TM},
\tau)\right\}^{(11)}=z_0(8\delta_2)^3+z_1 (8\delta_2)\epsilon_2,\ee
and by (3.32) and Theorem 4.1, \be
\left\{\mathrm{CS}\Phi_L(\nabla_0^{TM}, \nabla_1^{TM},
\tau)\right\}^{(11)}=2^6[z_0(8\delta_1)^3+z_1
(8\delta_1)\epsilon_1], \ee where, by comparing the
$q^{1/2}$-expansion coefficients in (4.22), \be z_0=-\left\{\int_0^1
\widehat{A}(TM,
\nabla_t^{TM})\mathrm{tr}\left[A\left(\frac{1}{2R_t^{TM}}-\frac{1}{8\pi\tan\left(
\frac{R_t^{TM}}{4\pi}\right)}\right)\right]dt\right\}^{(11)},\ee \be
\begin{split} z_1&=\left\{\int_0^1 \widehat{A}(TM,
\nabla_t^{TM})\mathrm{ch}\left(T_\CC M, \nabla_t^{T_\CC
M}\right)\mathrm{tr}\left[A\left(\frac{1}{2R_t^{TM}}-\frac{1}{8\pi\tan\left(
\frac{R_t^{TM}}{4\pi}\right)}\right)\right]dt\right.\\
&\left.+\int_0^1 \widehat{A}(TM,
\nabla_t^{TM})\mathrm{tr}\left[A\left(-\frac{1}{2\pi}\sin\left(\frac{R_t^{TM}}{2\pi}\right)+61\left(\frac{1}{2R_t^{TM}}-\frac{1}{8\pi\tan\left(
\frac{R_t^{TM}}{4\pi}\right)}\right) \right)
\right]dt\right\}^{(11)}.\end{split}\ee Plugging (4.24) and (4.25)
into (4.23) and comparing the constant terms of both sides, we
obtain that $$ \left\{\int_0^1 L(TM,
\nabla_t^{TM})\mathrm{tr}\left[A\left(\frac{1}{2R_t^{TM}}-\frac{1}{2\pi\sin\left(
\frac{R_t^{TM}}{\pi}\right)}\right)\right]dt\right\}^{(11)}=2^3(2^6z_0+z_1)
,$$ consequently \be \begin{split}&\left\{\int_0^1 L(TM,
\nabla_t^{TM})\mathrm{tr}\left[A\left(\frac{1}{2R_t^{TM}}-\frac{1}{2\pi\sin\left(
\frac{R_t^{TM}}{\pi}\right)}\right)\right]dt\right\}^{(11)}\\
=&\left\{\int_0^1 \widehat{A}(TM,
\nabla_t^{TM})\mathrm{ch}\left(T_\CC M, \nabla_t^{T_\CC
M}\right)\mathrm{tr}\left[A\left(\frac{1}{2R_t^{TM}}-\frac{1}{8\pi\tan\left(
\frac{R_t^{TM}}{4\pi}\right)}\right)\right]dt\right.\\
&\left.+\int_0^1 \widehat{A}(TM,
\nabla_t^{TM})\mathrm{tr}\left[A\left(-\frac{1}{2\pi}\sin\left(\frac{R_t^{TM}}{2\pi}\right)-3\left(\frac{1}{2R_t^{TM}}-\frac{1}{8\pi\tan\left(
\frac{R_t^{TM}}{4\pi}\right)}\right) \right)
\right]dt\right\}^{(11)}. \end{split}\ee We can view (4.26) as a
$11$-dimensional analogue of the miraculous cancellation formula
(3.38). We would like to remark that (4.26) can not be obtained
directly from (3.38) just by applying Chern-Simons transgression
to both sides as we pointed out right before Theorem 4.1.
Hopefully the odd dimensional cancellation formula (4.26) could
make sense in physics.

As for $\mathrm{CS}\Psi_W(\nabla_0^{TM}, \nabla_1^{TM}, \tau)$, we
have the following results.
\begin{theorem}Let $M$ be a $4k-1$ dimensional smooth flat manifold, $k \geq 2$, $\nabla_0^{TM}$ and $\nabla_1^{TM}$ be two flat
connections on $TM$, then for $2\leq i \leq k$,
$\left\{\mathrm{CS}\Psi_W(\nabla_0^{TM}, \nabla_1^{TM},
\tau)\right\}^{(4i-1)}$ is a modular form of weight 2i over $SL(2,
\mathbb{Z}).$
\end{theorem}
\begin{proof} Since both $\nabla_0^{TM}$ and $\nabla_1^{TM}$ are flat, we
have
$$(\nabla_0^{TM})^2=0,$$ $$(\nabla_1^{TM})^2=(\nabla_0^{TM}+A)^2=(\nabla_0^{TM})^2+[\nabla_0^{TM}, A]+A\wedge A=0, $$
which implies \be [\nabla_0^{TM}, A]=-A\wedge A.\ee One also has
$$R_t^{TM}=\left((1-t)\nabla_0^{TM}+t\nabla_1^{TM}\right)^2=(\nabla_0^{TM}+tA)^2=(\nabla_0^{TM})^2+t[\nabla_0^{TM}, A]+t^2A\wedge A.$$
Therefore \be R_t^{TM}=(t^2-t)A\wedge A.\ee Thus we obtain (cf. [30,
Lemma 1.7]) \be \mathrm{tr}[(R_t^{TM})^n]=(t^2-t)^n
\mathrm{tr}[A^{2n}]=\frac{(t^2-t)^n}{2}\mathrm{tr}[[A,A^{2n-1}]]=0,
\forall \,n \in \mathbf{Z}^+. \ee So it's not hard to see that
$$\mathrm{det}^{1\over
2}\left(f_{\Psi_W}\left(\frac{R_t^{TM}}{4\pi^2},
\tau\right)\right)=e^{{1\over
2}\mathrm{tr}\mathrm{ln}f_{\Psi_W}\left(\frac{R_t^{TM}}{4\pi^2},
\tau\right)}=1.$$ We therefore have
$$\mathrm{CS}\Psi_W(\nabla_0^{TM}, \nabla_1^{TM},
\tau)= {\frac{1} {8\pi^2}}\int_0^1
\mathrm{tr}\left[A\left(\frac{1}{\frac{R_t^{TM}}{4\pi^2}}-\frac{\theta'(\frac{R_t^{TM}}{4\pi^2},
\tau)}{\theta(\frac{R_t^{TM}}{4\pi^2}, \tau)}\right)\right] dt.$$

Then similar to (4.11), we have  \be \begin{split}
&\left\{\mathrm{CS}\Psi_W(\nabla_0^{TM},
\nabla_1^{TM}, -{1}/{\tau})\right\}^{(4i-1)}\\
=& \left\{{\frac{1} {8\pi^2}}\int_0^1
\mathrm{tr}\left[A\left(\frac{1}{\frac{R_t^{TM}}{4\pi^2}}-\frac{\theta'(\frac{R_t^{TM}}{4\pi^2},
-{1}/{\tau})}{\theta(\frac{R_t^{TM}}{4\pi^2},
-{1}/{\tau})}\right)\right] dt\right\}^{(4i-1)}\\
=& \left\{{\frac{1} {8\pi^2}}\int_0^1 \mathrm{tr}\left[A\tau
\left(\frac{1}{\frac{\tau
R_t^{TM}}{4\pi^2}}-\frac{\theta'(\frac{\tau R_t^{TM}}{4\pi^2},
\tau)}{\theta(\frac{\tau R_t^{TM}}{4\pi^2}, \tau)}\right)-A
\left(2\pi\sqrt{-1} \frac{\tau R_t^{TM}}{4\pi^2}\right)\right]
dt\right\}^{(4i-1)}\\
=&\left\{{\frac{1} {8\pi^2}}\int_0^1 \mathrm{tr}\left[A\tau
\left(\frac{1}{\frac{\tau
R_t^{TM}}{4\pi^2}}-\frac{\theta'(\frac{\tau R_t^{TM}}{4\pi^2},
\tau)}{\theta(\frac{\tau R_t^{TM}}{4\pi^2}, \tau)}\right)\right]
dt\right\}^{(4i-1)}, \end{split}\ee since $4i-1 \geq 7$ while
$\int_0^1 \mathrm{tr}\left[-A \left(2\pi\sqrt{-1} \frac{\tau
R_t^{TM}}{4\pi^2}\right)\right]$ is only a degree 3 form.

Thus we obtain that \be \begin{split}
&\left\{\mathrm{CS}\Psi_W(\nabla_0^{TM},
\nabla_1^{TM}, -{1}/{\tau})\right\}^{(4i-1)}\\
=&\left\{{\frac{1} {8\pi^2}}\int_0^1 \mathrm{tr}\left[A\tau
\left(\frac{1}{\frac{\tau
R_t^{TM}}{4\pi^2}}-\frac{\theta'(\frac{\tau R_t^{TM}}{4\pi^2},
\tau)}{\theta(\frac{\tau R_t^{TM}}{4\pi^2}, \tau)}\right)\right]
dt\right\}^{(4i-1)},\\
=&\tau^{2i}\left\{{\frac{1} {8\pi^2}}\int_0^1 \mathrm{tr}\left[A
\left(\frac{1}{\frac{ R_t^{TM}}{4\pi^2}}-\frac{\theta'(\frac{
R_t^{TM}}{4\pi^2}, \tau)}{\theta(\frac{ R_t^{TM}}{4\pi^2},
\tau)}\right)\right] dt\right\}^{(4i-1)}\\
=&\tau^{2i}\mathrm{CS}\Psi_W(\nabla_0^{TM}, \nabla_1^{TM}, \tau).
\end{split} \ee
It's also not hard to see that \be \mathrm{CS}\Psi_W(\nabla_0^{TM},
\nabla_1^{TM}, \tau+1)=\mathrm{CS}\Psi_W(\nabla_0^{TM},
\nabla_1^{TM}, \tau).\ee So $\left\{\mathrm{CS}\Psi_W(\nabla_0^{TM},
\nabla_1^{TM}, \tau)\right\}^{(4i-1)}$ is a modular form of weight
$2i$ over $SL(2, \mathbb{Z})$.
\end{proof}
\begin{remark}In Theorem 4.2, $i$ has to be greater than or
equal to 2 to get a modular form of weight 2i over $SL(2,
\mathbb{Z})$. This agrees with a known fact in number theory that
there is no nontrivial modular form of weight 2 over $SL(2,
\mathbb{Z})$.
\end{remark}

In particular, by Theorem 4.2,
$\left\{\mathrm{CS}\Psi_W(\nabla_0^{TM}, \nabla_1^{TM},
\tau)\right\}^{(7)}$ is a weight 4 modular form over $SL(2,
\mathbb{Z})$. Actually, \be
\begin{split} &\left\{\mathrm{CS}\Phi_L(\nabla_0^{TM}, \nabla_1^{TM},
\tau)\right\}^{(7)}\\
=& \left\{{\frac{1} {8\pi^2}}\int_0^1
\mathrm{tr}\left[A\left(\frac{1}{\frac{R_t^{TM}}{4\pi^2}}-\frac{\theta'(\frac{R_t^{TM}}{4\pi^2},
\tau)}{\theta(\frac{R_t^{TM}}{4\pi^2},
\tau)}\right)\right]dt\right\}^{(7)}\\
=& \left\{{\frac{1} {8\pi^2}}\int_0^1
\mathrm{tr}\left[A\left(\frac{1}{\frac{R_t^{TM}}{4\pi^2}}-\frac{\theta'(\frac{R_t^{TM}}{4\pi^2},
\tau)}{\theta(\frac{R_t^{TM}}{4\pi^2},
\tau)}\right)\right]dt\right\}^{(7)}\\
=&\frac{1} {512\pi^8}{1\over {3!}}\left.\frac{\partial^3}{\partial
z^3}\left(\frac{1}{z}-\frac{\theta'(z, \tau)}{\theta(z,
\tau)}\right)\right|_{z=0}\int_0^1\mathrm{tr}[A(R^{TM}_t)^3]dt\\
=&\frac{1} {512\pi^8}{1\over {3!}}\left.\frac{\partial^3}{\partial
z^3}\left(\frac{1}{z}-\frac{\theta'(z, \tau)}{\theta(z,
\tau)}\right)\right|_{z=0}\left(\int_0^1 (t^2-t)^3dt
\right)\mathrm{tr}[A^7].
\end{split}\ee
By direct computations, we can see that
$$ \frac{1} {512\pi^8}{1\over {3!}}\left.\frac{\partial^3}{\partial
z^3}\left(\frac{1}{z}-\frac{\theta'(z, \tau)}{\theta(z,
\tau)}\right)\right|_{z=0}\left(\int_0^1 (t^2-t)^3dt
\right)=-\frac{1}{3225600 \pi^4}+O(q).
$$
Let
$$ E_4(\tau)=1+240\sum_{n=1}^\infty\sigma_3(n)q^n,$$
be the Eisenstein series, which is a weight 4 modular form over
$SL(2, \mathbb{Z})$, where
$$\sigma_k(n)\triangleq \sum_{d|n}d^k.$$ It's a fact in number
theory that the space of weight 4 modular forms over $SL(2,
\mathbb{Z})$ has dimension 1. Thus
$$\left\{\mathrm{CS}\Phi_L(\nabla_0^{TM}, \nabla_1^{TM},
\tau)\right\}^{(7)}=-\frac{1}{3225600
\pi^4}E_4(\tau)\mathrm{tr}[A^7].$$

\section{Some discussions on loop space index theorem and flat vector bundles}
In this discussion section we first formally describe the loop space
version of the Atiyah-Singer index theory and it's relation with the
theory of elliptic genera following [21]. Then we compute the
Chern-Simons forms of two formal flat vector bundles on loop space.

Let $P$ be a principle $G$ bundle on $M$ and $E$ be an irreducible
positive energy representation of $\widetilde{L}G$, which is the
central extension of the loop group of $G$ [25]. We decompose $E$
according to the rotation action of the loop to get
$E=\sum_{\geq0}E_n$ where each $E_n$ is finite dimensional
representation of $G$. Constructing associated bundles to $P$ from
each $E_n$, which is still denoted by $E_n$, we get an element
$$\psi(P, E)=q^{m_\Lambda}\sum_n E_nq^n,$$ where $q=e^{2\pi \sqrt{-1}\tau}$ with $\tau \in
\mathbb{H}$ and $m_{\Lambda}$ is a rational number which is the so
called modular anomaly of the representation of $E$ ([17]).

For any positive integer $l$, the loop group
$\widetilde{L}\mathrm{Spin}(2l)$ has four irreducible level $1$
positive representations. Denote them by $S^+, S^-$ and $S_+, S_-$.
Let $\{\pm \alpha_j\}$ be the roots of $\mathrm{Spin}(2l)$. Then we
have the following normalized Kac-Weyl character formulas,
$$\chi_{S^+-S^-}=\prod_{j=1}^{l}\frac{\theta(\alpha_j,
\tau)}{\eta(\tau)}, \ \ \ \
\chi_{S^++S^-}=\prod_{j=1}^{l}\frac{\theta_1(\alpha_j,
\tau)}{\eta(\tau)},$$
$$\chi_{S_+-S_-}=\prod_{j=1}^{l}\frac{\theta_2(\alpha_j,
\tau)}{\eta(\tau)}, \ \ \ \
\chi_{S_++S_-}=\prod_{j=1}^{l}\frac{\theta_3(\alpha_j,
\tau)}{\eta(\tau)},$$ where
$\eta(\tau)=q^{1/{24}}\prod_{l=1}^\infty(1-q^l)$ is the Dedekind
eta-function [7].

Let $Q$ be the spin principle bundle associated to the tangent
bundle of $M$. We have the following,
$$\psi(Q, S^+-S^-)=q^{-k/12}(\Delta^+-\Delta^-)\otimes\bigotimes_{j=1}^{\infty}\Lambda_{-q^j}(TM), $$
$$\psi(Q, S^++S^-)=q^{-k/12}(\Delta^++\Delta^-)\otimes\bigotimes_{j=1}^{\infty}\Lambda_{q^j}(TM), $$
$$\psi(Q, S_+-S_-)=q^{-k/24}\otimes\bigotimes_{j=1}^{\infty}\Lambda_{-q^{j-{1/2}}}(TM), $$
$$\psi(Q, S_++S_-)=q^{-k/24}\otimes\bigotimes_{j=1}^{\infty}\Lambda_{q^{j-{1/2}}}(TM). $$
Thus we can view $S^+\pm S^-$ as the loop group analogues of the
finite dimensional spinor representations $\Delta^+\pm\Delta^-$,
where $\Delta^+$, $\Delta^-$ are the two irreducible spinor
representations of $\mathrm{Spin}(2l)$. Similar to the
$\widehat{A}$-class, the loop space $\widehat{A}$-class is defined
as
$$\widehat{\Theta}(M)=\frac{e(M)}{\mathrm{ch}(\psi(Q, S^+)-\psi(Q,
S^-))}=\eta(\tau)^k \cdot \prod_{j=1}^{k}\frac{x_j}{\theta(x_j,
\tau)}.$$ Therefore formally one defines the loop space Dirac
operator as
$$D^L=q^{-k \over {12}}D\otimes \bigotimes_{j=1}^{\infty}S_{q^j}(TM)$$
and the corresponding index formula
$$\mathrm{Ind}(D^L)=\int_M \widehat{\Theta}(M).$$ The twisted version
of this index theorem in this loop group setting is
$$\mathrm{Ind}(D^L\otimes \psi(P, E))=\int_M \widehat{\Theta}(M)\mathrm{ch}(\psi(P, E)).$$
In this sense we thus formally view $\psi(P, E)$ as a vector bundle
over the loop space $LM$.

To get modular forms instead of Jacobi forms, let's use virtual
versions of the above story. Denote by \be \mathcal{D}^L=D\otimes
\bigotimes_{j=1}^{\infty}S_{q^j}(\widetilde{T_\CC M}).\ee Let $V$ be
a $2l$ dimensional spin vector bundle over $M$. Physically,
$$\mathcal{V}:=\bigotimes_{j=1}^{\infty}\Lambda_{-q^{j-{1/2}}}(\widetilde{V_\CC}),\
\ \
\mathcal{V}':=\bigotimes_{j=1}^{\infty}\Lambda_{q^{j-{1/2}}}(\widetilde{V_\CC})$$
are virtual vector bundles over $LM$. Let $\nabla^{TM}$, $\nabla^V$
be two connections over $TM$, $V$ and $R^{TM}$, $R^V$ be their
curvatures respectively. By the Atiyah-Singer index theorem (also
cf. [20]), it's not hard to see that
$$ \mathrm{Ind}(\mathcal{D}^L\otimes \mathcal{V})=\int_M \mathrm{det}^{1\over
2}\left(\frac{R^{TM}}{4{\pi}^2}\frac{\theta'(0,\tau)}{\theta(\frac{R^{TM}}{4{\pi}^2},\tau)}
\frac{\theta_{2}(\frac{R^{V}}{4{\pi}^2},\tau)}{\theta_{2}(0,\tau)}\right)$$
and
$$\mathrm{Ind}(\mathcal{D}^L\otimes \mathcal{V}')=\int_M \mathrm{det}^{1\over
2}\left(\frac{R^{TM}}{4{\pi}^2}\frac{\theta'(0,\tau)}{\theta(\frac{R^{TM}}{4{\pi}^2},\tau)}
\frac{\theta_{3}(\frac{R^{V}}{4{\pi}^2},\tau)}{\theta_{3}(0,\tau)}\right).
$$
When $p_1(TM, \nabla^{TM})=p_1(V, \nabla^V)$,
$\mathrm{Ind}(\mathcal{D}^L\otimes \mathcal{V})$ and
$\mathrm{Ind}(\mathcal{D}^L\otimes \mathcal{V}')$ are modular forms
over $\Gamma^0(2)$ and $\Gamma_\theta$ respectively (cf. [21]).
Especially, taking $V=TM$ and $\nabla^V=\nabla^{TM}$, one has that
$\mathrm{Ind}(\mathcal{D}^L\otimes TLM)=\phi_W(M, \tau)$ and
$\mathrm{Ind}(\mathcal{D}^L\otimes TLM')=\phi_W'(M, \tau)$.

Now let's assume that $V$ is flat and $\nabla_0^V, \nabla_1^V$ are
two flat connections over $V$. Let $\nabla^{V}_i, i=0,1$ be two
connections on $V$ and $R^{V}_i, i=0,1$ be their curvatures
respectively. Let $\nabla_t^{V}=(1-t)\nabla_0^{V}+t\nabla_1^{V}$ and
$R_t^{V}$ be the corresponding curvature. Let $A=\nabla_1-\nabla_0
\in \Omega^1(M, \mathrm{End}(V)).$ One can lift these two
connections to $\mathcal{V}, \mathcal{V}'$ and denote them by
$\nabla_0^\mathcal{V}, \nabla_1^\mathcal{V}$ and
$\nabla_0^{\mathcal{V}'}, \nabla_1^{\mathcal{V}'}$. We heuristically
view $\mathcal{V}$ and $\mathcal{V}'$ as flat vector bundles on the
loop space $LM$.

Applying Theorem 2.2 similarly as in (4.1) to (4.4), one gets that
\be \begin{split} &\mathrm{ch}(\mathcal{V},
\nabla_1^\mathcal{V})-\mathrm{ch}(\mathcal{V},
\nabla_0^\mathcal{V})\\
=&\mathrm{det}^{1\over
2}\left(\frac{\theta_{2}(\frac{R_1^{V}}{4{\pi}^2},\tau)}{\theta_{2}(0,\tau)}\right)
-\mathrm{det}^{1\over
2}\left(\frac{\theta_{2}(\frac{R_0^{V}}{4{\pi}^2},\tau)}{\theta_{2}(0,\tau)}\right)\\
=&d\,\int_0^1 {\frac{1} {8\pi^2}}\mathrm{det}^{1\over
2}\left(\frac{\theta_{2}(\frac{R_t^{V}}{4{\pi}^2},\tau)}{\theta_{2}(0,\tau)}\right)\mathrm{tr}\left[A\left(\frac{\theta_2'(\frac{R_t^{V}}{4\pi^2},
\tau)}{\theta_2(\frac{R_t^{V}}{4\pi^2}, \tau)}\right)\right] dt.
\end{split}\ee
However by similar calculations as in the proof of Theorem 4.2,
one has
$$\mathrm{det}^{1\over
2}\left(\frac{\theta_{2}(\frac{R_t^{V}}{4{\pi}^2},\tau)}{\theta_{2}(0,\tau)}\right)=1.$$
Therefore we have \be \mathrm{ch}(\mathcal{V},
\nabla_1^\mathcal{V})-\mathrm{ch}(\mathcal{V},
\nabla_0^\mathcal{V})=d\,\mathrm{CS}(\mathcal{V},\nabla_0^\mathcal{V},\nabla_1^\mathcal{V},
\tau),\ee where \be
\mathrm{CS}(\mathcal{V},\nabla_0^\mathcal{V},\nabla_1^\mathcal{V},
\tau)\triangleq {\frac{1} {8\pi^2}}\int_0^1
\mathrm{tr}\left[A\left(\frac{\theta_2'(\frac{R_t^{V}}{4\pi^2},
\tau)}{\theta_2(\frac{R_t^{V}}{4\pi^2}, \tau)}\right)\right] dt \in
\Omega^{\mathrm{odd}}(M)[[q^{1\over 2}]].\ee Since $\nabla_0^V$ and
$\nabla_1^V$ are flat connections, $\mathrm{ch}(\mathcal{V},
\nabla_1^\mathcal{V})$ and $\mathrm{ch}(\mathcal{V},
\nabla_0^\mathcal{V})$ are both vanishing. Thus
$\mathrm{CS}(\mathcal{V},\nabla_0^\mathcal{V},\nabla_1^\mathcal{V},\tau)$
represents an element in $H^{\mathrm{odd}}(M, \mathbb{C})[[q^{1\over
2}]]$. One can similarly define \be
\mathrm{CS}(\mathcal{V}',\nabla_0^{\mathcal{V}'},\nabla_1^{\mathcal{V}'},
\tau)\triangleq {\frac{1} {8\pi^2}}\int_0^1
\mathrm{tr}\left[A\left(\frac{\theta_3'(\frac{R_t^{V}}{4\pi^2},
\tau)}{\theta_3(\frac{R_t^{V}}{4\pi^2}, \tau)}\right)\right] dt \in
\Omega^{\mathrm{odd}}(M)[[q^{1\over 2}]],\ee which also represents
an element in $H^{\mathrm{odd}}(M, \mathbb{C})[[q^{1\over 2}]]$.

Similarly as Theorem 4.1, we obtain that
\begin{theorem}Let $V$ be a $2l$ dimensional flat vector bundle over $M$ and $\nabla^{V}_0, \nabla^{V}_1$ be two flat connections on $V$,
then for any positive integer $i\geq 2$, we have
\newline 1)$\left\{\mathrm{CS}(\mathcal{V},\nabla_0^\mathcal{V},\nabla_1^\mathcal{V}, \tau)\right\}^{(4i-1)}$ is a modular form of weight 2i over
$\Gamma^0(2)$;
\newline $\left\{\mathrm{CS}(\mathcal{V}',\nabla_0^{\mathcal{V}'},\nabla_1^{\mathcal{V}'}, \tau)\right\}^{(4i-1)}$ is a modular form of weight 2i over
$\Gamma_\theta$;
\newline 2)\ the following equalities hold,
$$\mathrm{CS}(\mathcal{V},\nabla_0^\mathcal{V},\nabla_1^\mathcal{V}, \tau+1)=\mathrm{CS}(\mathcal{V}',\nabla_0^{\mathcal{V}'},\nabla_1^{\mathcal{V}'}, \tau).$$

\end{theorem}

Heuristically, (5.4) and (5.5) can be viewed as the Chern-Simons
transgressed forms of flat vector bundles over loop spaces. We
hope they could play some roles in the study of loop space vector
bundles.

\section{Acknowlegment}
The authors thank Professor Weiping Zhang for suggesting the
problems considered in this article. F. Han thanks Professor
Kefeng Liu and Professor Peter Teichner for many helpful
discussions and suggestions. Q. Chen thanks Professor Nicolai
Reshetikin for his interest and support. Part of our results were
announced on the Third Pacific Rim Conference held at Fudan
University in August 2005. The paper is finished when the second
author is visiting the Max-Planck-Institut f$\ddot{\mathrm{u}}$r
Mathematik at Bonn.

\end{document}